\documentclass[11pt]{article}

\usepackage{color}    
\usepackage{epsfig} 
\usepackage{epsf}            
\usepackage{float}           
\usepackage{amsthm}
\usepackage{graphics}
\usepackage{amsfonts}
\usepackage{amsxtra}
\usepackage{amsmath}
\usepackage{a4}
\usepackage{amssymb}
\usepackage{enumerate}

\newtheorem{theorem}{Theorem}
\newtheorem{lemma}[theorem]{Lemma}

\newtheorem{remark}[theorem]{Remark}

\numberwithin{theorem}{section}
\numberwithin{equation}{section}
\numberwithin{figure}{section}

\newcommand{\PP}{\mathbb{P}}
\newcommand{\ZZ}{\mathbb{Z}}

\newcommand{\EE}{\mathbb{E}\,}

\newcommand{\cP}{\mathcal{P}}

\newcommand{\cX}{\mathcal{X}}
\newcommand{\cU}{\mathcal{U}}

\newcommand{\tx}{{\tilde{x}}}

\allowdisplaybreaks

\begin{document}

\title{How to squeeze the toothpaste back into the tube}
\author{
\textbf{Pablo A.~Ferrari and James B.~Martin}
\\
\textit{Universidad de Buenos Aires and University of Oxford}
\\
\texttt{pferrari@dm.uba.ar, martin@stats.ox.ac.uk}
}
\date{1 March 2012}

\maketitle

\begin{abstract}
We consider \textit{bridges} for 
the simple exclusion process on $\ZZ$, either symmetric or asymmetric, 
in which particles jump to the right at rate $p$ and to the left at rate $1-p$. 
The initial state $O$ has all negative sites occupied and all non-negative sites empty. 
We study the probability that the process is again in state $O$ at time $t$, 
and the behaviour of the process on $[0,t]$ conditioned on being in state $O$ at time $t$. 
In the case $p=1/2$, we find that such a bridge typically goes a distance of order $t$ 
(in the sense of graph distance) from the initial state. 
For the asymmetric systems, $p\ne1/2$, we note an interesting duality which
shows that bridges with parameters $p$ and $1-p$ have the same distribution; the maximal 
distance of the process from the original state behaves like $c(p)\log t$ for some constant $c(p)$ 
depending on $p$. (For $p>1/2$, the front particle therefore travels
much less far than the bridge of the corresponding random walk,
even though in the unconditioned process the path of the front particle dominates a random walk.)
We mention various further questions.
\end{abstract}

\section{Introduction}\label{introductionsection}
How does the behaviour of a random process change if 
the process is conditioned to be back at its starting state at time $t$? 
Such questions were studied by Benjamini, Izkovsky and Kesten \cite{BIK}
for random walks on groups, and by 
Gantert and Peterson \cite{GantertPeterson} for 
one-dimensional random walks in a random environment. 

In this paper we consider the case of the simple exclusion process,
in both symmetric and asymmetric cases. Each site of $\ZZ$ may 
be occupied by particle or empty. Each particle 
attemps to jump one site to the right at rate $p$ and one
site to the left at rate $1-p$; the jump attempt succeeds if
the destination site is empty.

The time-0 configuration is the ``origin'' denoted $O$, in which
every negative site is occupied and every non-negative site is unoccupied.
In the interpretation of the process as a growth model in $\ZZ^2$
(``corner growth model''), $O$ corresponds
to a state where the positive quadrant of $\ZZ^2$ is unoccupied and the other three
quadrants are occupied.

We write $\PP_p$ or just $\PP$ for the law
of the process started from $O$, and
$\PP^x$ for the law of the process
started from a general state $x$. We write $X_t$ for the state of the process at time $t$.

The only states reachable from $O$ are those in which the number of
particles in $[0,\infty)$ is finite, and equal to the number of holes
in $(-\infty, 0)$. Write $\cU_0$ for the set of such states. For a state
$x\in \cU_0$, its \textit{distance} from the origin, denoted $D(x)$,  
is the minimum 
number of steps needed to get from $O$ to $x$ -- this is the distance
of $x$ from $O$ in the graph where two states are neighbours if 
the process can jump from one to the other.
Equivalently, this is the sum of the distance of each particle from its initial position. 
In the corner growth model interpreration, this is the area of the occupied
region in the positive quadrant.

We will study the probability that $X_t=O$, 
and the behaviour of process on $[0,t]$ conditioned on $X_t=O$
(that is, the behaviour of the \textit{bridge} of the process). 
The unconditioned process is transient for any $p\geq 1/2$. 

For $p=1/2$, we find lower and upper bounds of the form $e^{-c\sqrt{t}}$ for 
$\PP(X_t=O)$ (Theorem \ref{symmetricreturntheorem}). 
We show that the bridge typically reaches a maximum distance
of order $t$ from the origin during $[0,t]$ (Theorem \ref{symmetricdistancetheorem}).
This behaviour is on the same scale as that of the unconditioned process,
where $D_t/t\to 1/2$ as $t\to\infty$ (Lemma \ref{Dtlemma}). 

For the asymmetric case, we find an interesting duality between the cases with drift
to the left (recurrent, reversible in equilibrium
) and the cases with drift to the right (transient) -- conditioned
on $X_t=O$, the processes with parameter $p$ and $1-p$ have the same law (Theorem \ref{dualitytheorem}), 
and for $p>1/2$, the probability that $X_t=O$ grows like a constant
times $e^{-(2p-1)t}$. The maximum distance of the bridge from $O$ during $[0,t]$ is 
around $\left(\left|\log\frac{1-p}{p}\right|\right)^{-1}\log t$ 
from the origin during $[0,t]$; the same is true of the maximum distance of the 
rightmost particle from the origin (Theorem \ref{asymmetricdistancetheorem}). 
Note that the maximum deviation of the rightmost particle in the bridge 
has much smaller order ($\log t$) than the deviation of a random walk bridge 
with the same parameter (which would be $\sqrt{t}$) even though 
in the unconditioned process, the movement of the rightmost particle dominates 
a random walk. Loosely we could say that in the bridge, the rightmost particle
feels responsibility not only for returning itself to the origin at time $t$,
but also for preventing the other particles from wandering too far. 

In the next section, we introduce some notation and some give preliminary results. 
The results for the symmetric case are given in Section \ref{symmetricsection},
and those for the asymmetric case in Section \ref{asymmetricsection}.
Additional remarks and a variety of further questions are given in Section \ref{open}.

\section{Notation and Preliminaries}\label{notationsection}
Let $x\in\cU_0$. 
As above, let $D(x)$ be the graph distance of the state $x$ from $O$;
that is, the smallest number of jumps needed to get from $O$ to
$x$. Let $M(x)$ be the position of the furthest right particle in state $x$. Let $J(x)$ be
the number of leftward jumps available in state $x$, i.e.\ the
number of particles with a hole to their left. This is one fewer than
the number of rightward jumps available (the number of particles with a hole to their
right). Note that if $x\in\cU_0$, then certainly all of $D(x)$, $M(x)$ and $J(x)$
are finite. In particular, the rate of jumps out of state $x$ is finite,
and so the process is a Markov chain (i.e.\ has non-zero holding time in each state). 

Write $X_t$ for the state at time $t$. 
Write $D_t, M_t, J_t$ for $D(X_t), M(X_t), J(X_t)$. 

The following lemma collects together various useful facts:
\begin{lemma}\label{usefullemma}$\,$
\begin{enumerate}[(i)]
\item $J(x)\leq M(x)+1$. 
\item $J(x)\leq \sqrt{2D(x)}$.
\item The number of states $x$ with $D(x)=n$ is equal to the number of partitions of $n$, 
which is asymptotic to $\frac{\exp(\pi \sqrt{2/3} \sqrt{n})}{4n\sqrt{3}}$ as $n\to\infty$. 
\end{enumerate}
\end{lemma}
\begin{proof}
For $k\geq 1$, let $B^{(k)}(x)$ be the difference between the current position of the $k$th rightmost particle 
and its initial 
position (i.e.\ the minimal number of rightward jumps it must have made). So its current position is 
$B^{(k)}(x)-k$, since it started at site $-k$. Then $D(x)=\sum_k B^{(k)}(x)$.

Note that particle $k$ can make a leftward jump if and only if
$B^{(k)}>B^{(k+1)}$. But we have $B^{(1)}\geq B^{(2)}\geq\dots\geq 0$, so the total number of $k$ 
such that particle $k$ can make a leftward jump is at most $B^{(1)}$ which is equal to $M(x)+1$, giving (i). 

For (ii), note that if there are $j$ leftward jumps possible, 
then $B^{(k)}(x)$ takes at least $j$ different non-zero values as $k$ varies.
So $D(x)\geq 1+2+\dots+j\geq j^2/2$, which gives (ii).

Finally, note that $x$ is determined by the values $B^{(k)}(x), k\geq 1$, which form a decreasing sequence
of non-negative integers that sum to $D(x)$. This gives a bijection between states $x\in\cU_0$ 
and the set of integer partitions,
with the area (i.e.\ sum) of the partition equal to the distance $D(x)$, giving the correspondence in (iii).
The asymptotic enumeration of partitions of $n$ was first done by Hardy and Ramanujan in 1918.
\end{proof}

The following representation of $D(x)$ will be useful: 
\begin{lemma}\label{Dreplemma}
$D(x)=S^+(x) + S^-(x)$, where 
\begin{align*}
S^+(x)&=\sum_{i>0} i I(x\text{ has a particle at }i),\\
S^-(x)&=\sum_{j<0} (-j)I(x\text{ has a hole at }j).
\end{align*}
\end{lemma}
\begin{proof}
In the state $O$ we have $S^+=S^-=0$. 
Every jump right by a particle increases one of terms on the right-hand side by one,
and leaves the other unchanged. (If the particle jumps from $k$ to $k+1$ 
then the first sum increases if $k\geq 0$, and the second sum increases if $k<0$).
\end{proof}

We can represent the process $D_t$ in a natural way as the difference between two 
appropriately time-changed Poisson processes. 
Since the process $X_t$ is Markov, and jumps right (increasing the distance by 1) occur at rate $J(x)+1$
while jumps left (decreasing the distance by 1) occur at rate $J(x)$, we have the following result.
\begin{lemma}\label{Zreplemma}
We have the representation
\[
D_t=Z^+\left(\int_0^t [J_s + 1]ds\right) - Z^-\left(\int_0^t J_s ds\right)
\]
where $Z^+$ and $Z^-$ are independent Poisson processes with rates $p$ and $1-p$.
\end{lemma}

Finally we give two useful results concerning stochastic ordering of states and related coupling properties.
For two states $x$ and $\tx$ in $\cU_0$, we will write $\tx\preccurlyeq x$ and say that $x$ dominates $\tx$
if, for any site $i$, the number of particles to the right of $i$ is at least as large in $x$ as it is in $\tx$.
Equivalently, for each $k$, 
the $k$th rightmost particle in $\tx$ is no further right than the $k$th rightmost particle in $x$. 
Note that $O\preccurlyeq x$ for all $x\in\cU_0$.
\begin{lemma}\label{couplinglemma}
Let $x$ and $\tx$ be two states in $\cU_0$ such that $\tx\preccurlyeq x$. 
Then $\PP^{x}(X_t=O) \leq \PP^{\tx}(X_t=O)$.
\end{lemma}
\begin{proof}
Since $O$ is the minimal state, the result follows as soon as we have a coupling of the two processes which
preserves the domination. There are a variety of natural couplings which do the job; for example,
the well-known ``basic coupling'': each site $i$ carries a Poisson process $\cP^+_i$ of rate $p$ and a 
Poisson process $\cP^-_i$ of rate $1-p$, all of which are independent. We use the same collection of Poisson
processes to run both processes: whenever a point of $\cP^+_i$ (respectively $\cP^-_i$) occurs,
a jump is attempted from site $i$ to site $i+1$ (respectively site $i-1$), which succeeds whenever site $i$
is occupied and the destination site is unoccupied. With the jump attempts coupled in this way,
it is impossible for the 
$k$th particle in one process to overtake the $k$th particle in the other by a rightward jump, unless
the $(k-1)$st particle in that process has already overtaken the $(k-1)$st particle in the other process.
A similar observation holds for leftward jumps.
Hence if the 
initial states are ordered by $\preccurlyeq$, then this ordering is maintained for all times, as required.
\end{proof}

\begin{lemma}\label{dominationlemma}
The process conditioned on $X_t=O$ is dominated on the time interval $[0,t]$ by the unconditioned process. 
That is, there is a coupling of the two processes 
such that at all times in $[0,t]$, the state of the unconditioned
process dominates that of the conditioned process.
\end{lemma}
\begin{proof}
The conditioned process is itself a Markov chain (though of course not time-homogeneous). 
Suppose that the unconditioned process has rate $r$ of jumping from state $x$ to state $y$. 
From Bayes' rule, one sees that if the 
conditioned process is in state $x$ at time $s$, and $\PP(X_t=O | X_s=x) < \PP(X_t=O | X_s=y)$,
then the conditioned process has instantaneous rate greater than $r$ of jumping from state $x$ to state $y$.
Similarly if $\PP(X_t=O | X_s=y) < \PP(X_t=O | X_s=x)$ then the conditioned process will 
have instantaneous rate less than $r$ of jumping from $x$ to $y$.  

Now observe that Lemma \ref{couplinglemma} tells us that any rightward
jump reduces the probability of ending up in state $O$ at time $t$,
while any leftward jump increases that probability.  Hence in the
conditioned process, every possible rightward jump is occurring with
some instantaneous rate less than $p$, while every possible leftward
jump is occurring with some instantaneous rate greater than $1-p$.

Now we couple the two processes, using for example the basic coupling
as described in the proof of Lemma \ref{couplinglemma}. But in
addition we must now suppress some of the rightward jump attempts in the
conditioned process, and also add some extra leftward jump attempts in
that process.  Since the two initial states are the same, this
coupling will maintain domination of the state of the conditioned
process by that of the unconditioned process at all times.
\end{proof}

\section{$p=1/2$}\label{symmetricsection}

We will show upper and lower bounds of the form $e^{-c\sqrt{t}}$
for the probability of being at the origin at time $t$. We'll show that 
the bridge on $[0,t]$ typically reaches a maximal distance on order $t$ 
from the initial state.

\begin{lemma}\label{squaresumlemma}
If $p=1/2$ then
\[
\PP(X_t=O)=\sum_x \PP(X_{t/2}=x)^2.
\]
\end{lemma}

\begin{lemma}\label{Dtlemma}
If $p=1/2$ then $\frac{D_t}{t}\to\frac12$ a.s.\ as $t\to\infty$.
\end{lemma}

\begin{proof}
Arratia \cite{arratia} showed that 
\begin{equation}\label{arratia}
t^{-1/2}M_t - \sqrt{\log t} \to 0 \text{ a.s.\ as } t\to\infty.
\end{equation}
This is a much more precise estimate than we need. Since $J_t\leq M_t+1$ (Lemma \ref{usefullemma}(i)), 
we get $J_t/t^{3/4}\to 0$ a.s., 
and so 
\[
\frac{\int_0^t J_s ds}{t^{7/4}}\to 0 \text{ a.s.}
\]
Recall the processes $Z^+$ and $Z^-$ from Lemma \ref{Zreplemma}, 
which are now both Poisson processes of rate $1/2$.
Of course, they are not independent of $(J_t)$, but 
we can apply, for example, the law of the iterated logarithm 
to give 
\begin{align*}
\frac{
Z^+\left(\int_0^t J_s ds + t\right)-\frac12\left(\int_0^t J_s ds + t\right)
}{t}
&\to 0 \text{ a.s.}\\
\frac{
Z^-\left(\int_0^t J_s ds \right)-\frac12\left(\int_0^t J_s ds \right)
}{t}
&\to 0 \text{ a.s.}
\end{align*}
So indeed 
\begin{align*}
\frac{D(t)}{t}&=\frac{Z^+\left(\int_0^t J_s ds + t\right) - Z^-\left(\int_0^t J_s ds\right)}{t}\\
&\to \frac12 \text{ a.s.},
\end{align*}
as required.
\end{proof}

\begin{lemma}\label{typicalsetlemma}
Let $p=1/2$. There are sets $\cX_t$ with the following properties:
\begin{itemize}
\item[(i)] $\PP(X_t\in\cX_t)\to 1$ as $t\to\infty$;
\item[(ii)] 
\[
|\cX_t|\approx \exp\left(\frac{\pi}{\sqrt{3}}\sqrt{t}\right),\] 
in the sense that
\[
\lim_{t\to\infty} \frac{\log(|\cX_t|)}{\sqrt{t}}=\frac{\pi}{\sqrt{3}}.
\]
\end{itemize}
\end{lemma}

\begin{proof}
This follows from the facts that $D_t/t$ converges a.s.\ to $1/2$ (Lemma \ref{Dtlemma}), and 
that the number of states at distance $n$ is asymptotic to $\exp(\pi\sqrt{2/3}\sqrt{n})$,
up to polynomial corrections (Lemma \ref{usefullemma}(iii)). 
\end{proof}


\begin{theorem}\label{symmetricreturntheorem}
Let $p=1/2$.
There exist $c, c'>0$ such that
\begin{equation}
\label{upperlower}
-c\leq \liminf_{t\to\infty} \frac{\log \PP(X_t=O)}{\sqrt{t}} 
\leq \limsup_{t\to\infty} \frac{\log \PP(X_t=O)}{\sqrt{t}}
\leq -c'.
\end{equation}
In fact one can take $c=\pi/\sqrt{6}\approx 1.2825$ and 
$c'=\int_0^\infty -\log\Phi(x)dx\approx 0.4775$.
\end{theorem}

\begin{proof}
Consider the sets $\cX_t$ from Lemma \ref{typicalsetlemma}. 
Using Lemma \ref{squaresumlemma}, we have
\begin{align*}
\PP(X_t=O)&=\sum_x \PP(X_{t/2}=x)^2\\
&\geq \sum_{x\in|\cX_{t/2}|} \PP(X_{t/2}=x)^2\\
&\geq \frac{\PP(X_{t/2}\in\cX_{t/2})}{|\cX_{t/2}|}.
\end{align*}
Then from the properties of the set $\cX_t$ in Lemma \ref{typicalsetlemma}, we get  
\begin{equation}\label{below}
\liminf_{t\to\infty} \log\PP(X_t=O)/\sqrt{t} \geq -\pi/\sqrt{6}.
\end{equation}

In the other direction, we use the negative dependence properties of the 
symmetric exclusion process. First consider the ``stirring'' representation
of the process. Each particle can be seen as performing a simple symmetric random walk 
on $\ZZ$ --- for this interpretation, we imagine that each pair of neighbouring particles
exchange their positions at rate $1/2$, as well as the normal movement of particles 
into empty spaces. Of course these random walks are not independent; rather they are negatively
dependent, including in the following sense. 

Consider any initial configuration which includes particles at the sites $y_1, y_2,\dots, y_k$,
and let $A$ be any subset of $\ZZ$. Let $P_1$ be the probability that all of the $k$ particles
have positions in the set $A$ at time $t$. 

Consider also a collection of independent simple symmetric random walks, starting at positions
$y_1, y_2, \dots, y_k$. Let $P_2$ be the probability that all of these particles have positions
in the set $A$ at time $t$. 

Then $P_1\leq P_2$. This follows, for example, from a more general result such as Proposition
1.7 in Chapter VIII of \cite{Liggettbook}. (That result applies to 
expectations of more general positive definite functions of the positions of the $k$ particles; 
here we consider simply the function which is the product of $k$ copies of the identity function 
of the set $A$).

Hence an upper bound $P(X_t=O)$ is given by the probability that independent simple symmetric 
random walks started at
positions $-1, -2, \dots, -k$ are all at negative positions at time $t$. Write $W^{(y)}(t)$ 
for the position at time $t$ of such a walk started at position $y$. 
Since $k$ can be
taken arbitrarily large, we have
\begin{align*}
\PP(X_t=O)
&\leq \prod_{y=-1,-2,\dots} \PP\left(W_t^{(y)}<0\right)\\
&=\prod_{v=1,2,\dots} \PP\left(W_t^{(0)}<v\right).
\end{align*}
For any fixed $K$, the central limit theorem gives 
$\PP\left(W^{(0)}_t< v\right)=\Phi\left(\frac{v}{\sqrt{t}}\right)+o(1)$ 
uniformly in $v\in[1,K\sqrt{t}]$ as $t\to\infty$, so that
\begin{align*}
\limsup_{t\to\infty}\frac{\log\PP(X_t=O)}{\sqrt{t}}
&\leq \limsup_{t\to\infty} \sum_{v=1}^{\lfloor K\sqrt{t}\rfloor} 
\frac{\log\Phi\left(\frac{v}{\sqrt{t}}\right)}{\sqrt{t}}\\
&=\int_0^K\log\Phi(x) dx.
\end{align*}
Taking $K$ arbitrarily large gives the upper bound in 
(\ref{upperlower}) with $c'=\int_0^\infty -\log\Phi(x) dx$, as required. 
\end{proof}

\begin{theorem}\label{symmetricdistancetheorem} Let $p=1/2$.
If $\alpha<\frac{3}{4\pi^2}$ and $\beta>\frac14$, then 
\begin{equation}\label{symmetricdistance}
\PP\big(\alpha t<\max_{s\in[0,t]} D_s<\beta t | X_t=O\big)\to 1 \text{ as } t\to\infty.
\end{equation}
\end{theorem}

\begin{proof}
For the lower bound, we'll show that if $\alpha<\frac{3}{4\pi^2}$,   then
\begin{equation}\label{compare}
\PP(X_t=O, D_s<\alpha t \text{ for all } 0<s<t) = o(\PP(X_t=O).
\end{equation}
We look for an estimate for the left hand side of (\ref{compare}) to compare with the 
lower bound for $\PP(X_t=O)$ in (\ref{upperlower}). 

The idea is as follows: because there is always one more jump to the
right than jump to the left available, $D_t$ has mean $t/2$. So to see
$X_t=O$, i.e.\ $D_t=0$, we need a significant deviation from the mean in
either the process of rightward jumps or of leftward jumps.  It is hard for this 
to happen unless the expected number of such jumps is rather
large; this makes it unlikely for $D(s)$ to remain small for the
whole interval, since when $D(s)$ is small there are few jumps
available.

Using Lemma \ref{usefullemma}(ii), we have
\begin{align*}
\PP\left(X_t=O, \max_{s\in[0,t]}D_s\leq\alpha t\right)
&\leq
\PP(D_t=0, J_s\leq\sqrt{2\alpha t} \text{ for all } 0<s<t)\\
&\leq \PP(\text{ for some } u<\sqrt{2\alpha}t^{3/2}, \,\, Z^+(u+t)-Z^-(u)=0),
\end{align*}
where $Z^+$ and $Z^-$ are independent Poisson processes of rate $1/2$, 
as in Lemma (\ref{Zreplemma}). 
Then by the result in Lemma \ref{Zlemma} below,
\begin{equation}\label{upperbound}
\limsup_{t\to\infty} \frac{\log\PP\left(X_t=O, \max_{s\in[0,t]}D_s<\alpha t\right)}{\sqrt{t}}
\leq 
-\frac{1}{2\sqrt{2\alpha}}.
\end{equation}
Comparing the exponents in (\ref{upperlower}) and (\ref{upperbound}), we see
that (\ref{compare}) holds whenever 
$\frac{\pi}{\sqrt{6}}<\frac{1}{2\sqrt{2\alpha}}$, i.e.\ $\alpha<\frac{3}{4\pi^2}$, as required.

For the upper bound, we'll first show that
\begin{equation}
\label{firsthalf}
\PP(D_s< \beta t \text{ for all } 0<s<t/2 | X_t=O) \to 1\text{ as } t\to\infty.
\end{equation}
From the domination result in Lemma \ref{dominationlemma}, it is enough to obtain the equivalent estimate for 
the unconditioned process, i.e.\ that 
\begin{equation}
\label{newfirsthalf}
\PP(D_s< \beta t \text{ for all } 0<s<t/2) \to 1\text{ as } t\to\infty.
\end{equation}
But since $\beta>1/4$, this follows from the fact that $D_t/t \to 1/2$ a.s.\ as given by Lemma \ref{Dtlemma}.

Finally observe that since we start in state $O$ and condition on finishing in state $O$, 
the distribution of the path on $[0,t]$ is invariant under time-reversal. Hence from 
(\ref{firsthalf}) we also get immediately that
\begin{equation}
\label{secondhalf}
\PP(D_s< \beta t \text{ for all } t/2<s<t | X_t=O) \to 1\text{ as } t\to\infty.
\end{equation}
The convergence in (\ref{symmetricdistance}) now follows by combining (\ref{compare}), (\ref{firsthalf})
and (\ref{secondhalf}).
\end{proof}

\begin{lemma}\label{Zlemma}
Let $Z^+$ and $Z^-$ be independent Poisson processes, both of rate $1/2$. Let $\gamma>0$. 
Then 
\[
\limsup_{t\to\infty} 
\frac{
\log\PP\left(Z^+(u+t)-Z^-(u)\leq 0 \text{ for some } u<\gamma t^{3/2}\right)
}{t^{1/2}}
\geq
-\frac{1}{2\gamma}.
\]
\end{lemma}

\begin{proof}
Write $T=\gamma t^{3/2}$. Fix any $\delta>0$. 
We will consider the process $Y(u)=Z^+(u+t)-Z^-(u)$, $u\geq 0$, and also its jump chain.

Consider the following three events:
\begin{itemize}
\item[(a)] The initial value $Y(0)=Z^+(t)$ is less than $(1-\delta)t/2$.
\item[(b)] The process $Y$ makes more than $(1+\delta)T$ jumps in $[0,T]$.
\item[(c)] The process $Y$ goes at least $(1-\delta)t/2$ below its initial value
at some time during its first $(1+\delta)T$ jumps.
\end{itemize}

If none of these events occur, then $Z^+(u+t)-Z^-(u)>0$ for all $u\leq T$. 

Since $Z^+$ and $Z^-$ are both Poisson processes of rate $1/2$, the event in (a)
decays exponentially in $t$, and the event in (b) decays exponentially in $T$ and hence
faster than exponentially in $t$.

It remains to estimate the probability of the event in (c). The jump
chain is a martingale, with jumps of $\pm 1$. We consider the jump
chain stopped at the first time it goes at least $(1-\delta)t/2$ below
its initial value (the stopped chain is also a martingale with jumps
bounded by 1).  Applying the Azuma-Hoeffding inequality to that chain,
we can bound the probability of the event in (c) above by
\[
\exp\left(-\frac{\big[(1-\delta)t/2]^2}{2(1+\delta)T}\right).
\]

Combining the estimates for (a), (b) and (c), using $T=\gamma
t^{3/2}$, and noting that $\delta$ can be arbitrarily small, we obtain
the upper bound in the statement of the lemma.
\end{proof}

\section{$p\ne 1/2$}\label{asymmetricsection}
\subsection{Duality}
Consider a path $x_{[0,t]}=\{x_s, 0\leq s\leq t\}$ of the process on the time
interval $[0,t]$.  
The number of rightward jumps available at time $s$ is $J(x_s)+1$ and the number 
of leftward jumps available is $J(x_s)$. Let $R$ be the total number of 
rightward jumps made and $L$ the total number of leftward jumps made. 
In the representation of Lemma \ref{Zreplemma},
the path corresponds to a realisation with $R$ points in the interval
$\left[0,\int_0^t\left(J(x_s)+1\right)ds\right]$ in the process $Z_+$,
and $L$ points in the interval $\left[0,\int_0^t J(x_s)ds\right]$ in the process $Z_-$.
Then the likelihood of the path is
\begin{align*}
w(x_{[0,t]})
&=
p^R (1-p)^L 
\exp\left(-p\int_0^t \left[J(x_s)+1\right]ds\right)
\exp\left(-(1-p)\int_0^t J(x_s)ds\right)
\\
&=
p^R (1-p)^L 
e^{-pt}
\exp\left(-\int_0^t J(x_s)ds\right).
\end{align*}
If the path begins and ends at $O$, then $R=L$.
Then exchanging $p$ and $1-p$ simply multiplies the likelihood of the path by $e^{(2p-1)t}$.
Since, for given $t$, this scaling is the same for all paths that begin and end at $O$, 
we can integrate over all such paths to give:
\begin{theorem}\label{dualitytheorem}
$\,$
\begin{itemize}
\item[(a)]
$\PP_p(X_t=O) = e^{-(2p-1)t}\PP_{1-p}(X_t=O)$.
\item[(b)]
Conditioned on the event $X_t=O$, the processes with parameters $p$ and $1-p$ 
have identical laws on the interval $[0,t]$.
\end{itemize}
\end{theorem}

\begin{remark}
The duality above can be described as follows: the distribution of a bridge
from $O$ to $O$ on the interval $[0,t]$ is the same for the process with
parameter $p$ as it is for the process with parameter $1-p$.

We can compare this to well-known results about simple random walk.
Consider instead a simple random walk on $\ZZ_+$, with jumps to the right at rate $p$, 
and with jumps to the left at rate $1-p$ except at the state 0. In this case,
the distribution on excursions from $O$ to $O$ on $[0,t]$ is the same if $p$ is replaced
by $1-p$; in fact, it's the same for all $p\in(0,1)$. However, the distribution
on \textit{bridges} is not the same. For $p<1/2$, the unconditioned process is recurrent
and the bridge tends to stay close to the origin. For $p>1/2$, the bridge feels a penalty
for spending time at the origin and so tends not to visit the origin except
near the beginning and end of the interval.
(If there is no barrier at 0, so the walk jumps left at rate $1-p$ and right
at rate $p$ from all sites, then the distributions on both bridges and excursuions
don't depend on $p$).
\end{remark}

\subsection{Stationarity}
For any $p$ and any $\lambda\in(0,1)$, 
the product measure under which each site is occupied independently with
probability $\lambda$ is a stationary distribution for the process. 
(This stationary distribution corresponds to a reversible process only if $p=1/2$).

When $p\ne 1/2$, there are further stationary distributions known as 
\textit{blocking measures} (see e.g.\ \cite{BraLigMou}). 
Consider the product measure $\mu_p$ under which site $i$ is occupied 
with probability
\[
\frac{1}{1+\left(\frac{1-p}p\right)^i}.
\]
This is stationary for the process. 
Note that if $p<1/2$, then this distribution is concentrated
on states where $N_+$, the number of occupied non-negative sites,
and $N_-$, the number of unoccupied positive sites, are both finite. 
However, the state space is not irreducible; the quantity $N_+-N_-$
is preserved by the dynamics. Since we are interested in the process
started from the state $O$, we look at the process restricted
to the space $\cU_0=\{N_+-N_-=0\}$. This process is now irreducible;
the process on $\cU_0$ is a positive recurrent Markov chain,
with stationary distribution $\mu^{\cU_0}_p$ which is simply the distribution
$\mu_p$ conditioned on $N_+-N_-=0$ (which is an event of positive probability under $\mu_p$).





\begin{remark}
We mention one further interpretation of the duality above. 
Let $p<1/2$. Let $\alpha=\mu_{p}^{\cU_0}(O)=\mu_p(O)/\mu_p(\cU_0)$.

Since the process with parameter $p$ is irreducible and positive recurrent, we have that as $t\to\infty$, 
$\PP_p(X_t=O)\to\alpha$. We can then restate the first part of Theorem \ref{dualitytheorem} to 
describe the behaviour of the system with the drift to the right:
\[
\PP_{1-p}(X_t=0)\sim e^{-(1-2p)t}\alpha \text{ as } t\to\infty.
\]
Consider any state $x$. Any path from $x$ to $O$ has $D(x)$ more leftward steps
than it has rightward steps. Hence by the same argument that led to 
Theorem \ref{dualitytheorem}, 
writing $\PP_{1-p}^x$ for the process started from state $x$,
\[
\PP^x_{1-p}(X_t=O)\sim e^{-(1-2p)t}\alpha\left(\frac{p}{1-p}\right)^{D(x)}
\text{ as } t\to\infty.
\]
In that sense we can see $\left(\frac{p}{1-p}\right)^{D(x)}$
as the $h$-transform corresponding to conditioning 
on the event of returning to the origin at some distant time $t$. 
\end{remark}

\subsection{Distance from the origin}
We now look at how far the process moves from the origin 
(in the sense of the distance function $D$)
when conditioned to return to $O$ at time $t$. 

\begin{theorem}\label{asymmetricdistancetheorem}
Let $p\ne 1/2$. Let $c(p)=\frac{1}{\left|\log\frac{1-p}p\right|}$.
If $\epsilon>0$, then
as $t\to\infty$,
\begin{equation}\label{asymmetricdistance}
\PP_p\left(
(1-\epsilon)c(p)\log t
\leq \max_{s\in(0,t)} M_s \leq \max_{s\in(0,t)} D_s 
\leq 
(1+\epsilon){c(p)}\log t 
\,\big|\,X_t=O\right)
\to 1.
\end{equation}
\end{theorem}
The rest of this section contains the proof of Theorem \ref{asymmetricdistancetheorem}.
In the light of Theorem \ref{dualitytheorem}, 
it will be enough to show the result for the case $p<1/2$. 

Since the process is positive recurrent,
note that the event $X_t=O$ occurs with uniformly
positive probability under $\PP_p$. Hence if we have any collection of events 
$A_t$ such that $\PP_p(A_t)\to 1$ as $t\to\infty$,
then also $\PP_p(A_t|X_t=O)\to 1$. Hence it will be enough to show that
the event in (\ref{asymmetricdistance})
has unconditioned probability tending to 1 under $\PP_p$. 

For the lower bound, note that the position of the front particle is bounded below
by a random walk reflecting at 0, which jumps rightwards from any site at rate $1-p$
and leftwards from any site except 0 at rate $p$. 
Such a walk has stationary distribution 
\[
\pi(i)=\frac{1-2p}{1-p}\left(\frac{p}{1-p}\right)^i, 
\, i=0,1,2,\dots.
\]
From any state $i>0$, consider the probability of hitting site $i$ again before next hitting 0. 
This is no larger than the probability, in the walk without reflecting boundary at 0, 
of never hitting $i$ again (i.e.\ of hitting $-\infty$ before $i$). Hence in particular
this is bounded above uniformly for all sites $i$. 

As a result, the expected time to hit site $i$ from 0 is bounded above
by a constant times $1/\pi(i)$. 

Let $i(t)$ be a site which depends on $t$. If
$i(t)\leq O({1-\epsilon})c(p)\log t$ for all $t$, 
then $1/\pi(i(t))=o(t)$ as $t\to\infty$, and so the expected time to hit site $i(t)$ from 0 is $o(t)$. 
Thus the probability that the random walk hits site $i(t)$ in the time interval $[0,t]$ goes to 1 as $t\to\infty$.
The same must be true of the front particle in the exclusion process. 
So indeed the first inequality in (\ref{asymmetricdistance}) holds with probability tending to 1 as $t\to\infty$.

\newcommand{\Ppstat}{{\PP_p^{\text{stat}}}}

We turn to the upper bound. Recall from Lemma \ref{dominationlemma}
that the conditioned process with law $\PP_p(.|X_t=O)$ is dominated
by the unconditioned process with law $\PP_p(.)$. 
This is turn is dominated by the process started in stationarity (i.e.\ with $X_0$ distributed according to $\mu_p^{\cU_0}$),
whose law we denote by $\Ppstat$.
Hence it will be enough to show that the upper inequality in (\ref{asymmetricdistance}) holds with 
probability tending to 1 under $\Ppstat$.
 
We look at the tail of $D(x)$ for a state $x$ drawn from the stationary distribution $\mu_p^{\cU_0}$:

\begin{lemma}\label{exptaillemma}
Suppose $c>c(p)$. Then 
  $\mu_p^{\cU_0}\left(x:D(x)>c\log t\right)=o(1/t)$ as $t\to\infty$.
\end{lemma}
\begin{proof}
We know that $\mu_p^{\cU_0}(.)=\mu_p(.)/\mu_p(\cU_0)$,
so it will be enough to show the equivalent bound under $\mu_p$ rather than $\mu_p^{\cU_0}$.
Write $\EE_p$ for the expectation under $\mu_p^{\cU_0}$. Then it will be enough to show 
that $\EE_p e^{D(x)/c}<\infty$. 

Lemma \ref{Dreplemma} gives $D=S^+ + S^-$, where
$S^+=\sum_{i>0}iI(\text{site $i$ occupied})$
and $S^-=\sum_{j<0}(-j)I(\text{site $j$ unoccupied})$. Under $\mu_p$, $S^+$ and $S^-$ are independent. 
So it will be enough to show that $\EE_p e^{S^+/c}<\infty$ and $\EE_p e^{S^-/c}<\infty$.
The arguments for $S^+$ and $S^-$ are essentially identical; we write the one for $S^+$. 

Under $\mu_p$, the sites are occupied independently, and site $i$ is occupied with 
probability $r_i=\left(1+\left(\frac{1-p}{p}\right)^i\right)<\left(\frac{p}{1-p}\right)^{-i}$.
So 
\begin{align*}
\EE_p e^{S^+/c}
&=\prod_{i>0}\left(1-r_i+r_i e^{i/c}\right)\\
&\leq\prod_{i>0}\left(1+\left(\frac{pe^{1/c}}{1-p}\right)^i\right).
\end{align*}
Since $e^{1/c}<\frac{1-p}{p}$, the right hand side is finite, as required.
\end{proof}

\begin{lemma}
For any $\epsilon>0$,
\[
\Ppstat\left(
\max_{s\in(0,t)}D(s) > (1+\epsilon)c(p)\log t
\right)
\to 0 \text{ as } t\to\infty.
\]
\end{lemma}
\begin{proof}
Consider any distance $d$ and any $\delta>0$.
Suppose the process is in state $x$ with $D(x)=d$. 
Now consider the amount of time until the process first
reaches any state $x'$ with $D(x')\leq (1-\delta)d$. 
For this to happen, we must see at least $\delta d$ leftward 
jumps from states $y$ with $D(y)\leq d$. 
But if $D(y)\leq d$ then $J(y)$,
the number of leftward jumps available from state $y$,
is at most $\sqrt{2d}$. 
So the time until we see at least $\delta d$ leftward jumps
from states $y$ with $D(y)\leq d$ 
dominates the time to see $\delta d$ points in a Poisson process of rate $(1-p)\sqrt{2d}$. 
As $d\to\infty$, the probability that this time is less than 2 goes to 0;
certainly for large enough $d$, if $D_s=d$ then with probability 
at least $1/2$, $D_{s+u}\geq (1-\delta)d$ for all $u\in[0,2]$.

Hence, for any large enough $d$ and any $t$, 
the expected amount of time in $[0,t+2]$ that
the process spends in states at distance $(1-\delta)d$ is at least
as big as the probability that the process visits a state at distance $d$
during $[0,t]$. 

Now given $\epsilon>0$, choose $\delta>0$ with $(1-\delta)(1+\epsilon)=1+\epsilon/2$.
From Lemma \ref{exptaillemma} we see that in stationarity, 
the expected amount of time spent in states at distance 
$(1+\epsilon/2)c(p)\log t$
in $[0,t+2]$ is $o(1)$ as $t\to\infty$. 
Therefore the probability of visiting a state at distance 
${(1+\epsilon)}c(p)\log t$
during $[0,t]$ is also $o(1)$.
\end{proof}

This completes the proof of Theorem \ref{asymmetricdistancetheorem}.

\section{Further remarks and questions}\label{open}
\begin{enumerate}[(i)]
\item
Note that when $p\ne 1/2$, the maximum distance of the bridge from the origin
is asymptotically no larger than the maximum position of the front particle alone. 
So typically at some point one sees the front particle reach distance around $c(p)\log t$,
while all the other particles are at positions $o(\log t)$. 
One could also ask about the maximal distance from the origin
reached by the other particles. By similar methods one can show that with high
probability the second particle reaches distance approximately $(c(p)/2)\log t$
(at which time again the distance of the whole state from the origin is around $c(p)\log t$
and all the particles except the front two are at positions $o(\log t)$). 
Similarly the third particle will reach distance approximately $(c(p)/3)\log t$ at some
point, and so on. In fact, any configuration with $D(x)<(1-\epsilon)c(p)\log t$ has 
expected hitting time $o(t)$ from O, and so will typically be visited during the time interval
$(0,t)$ (both in the case of the bridge and in the case of the unconditioned process).
\item
From the duality with the recurrent case $p<1/2$, it is straightforward
to sample from the distribution of the bridge in the case $p>1/2$. 
How would one sample from the distribution in the symmetric case?
\item
What can be said about the maximal distance of the rightmost particle from the origin
in the case $p=1/2$? By symmetry between holes and particles, the process $M_s$ has
the same distribution as the process $H_s$, where $-(H_s+1)$ is the furthest left
position of a hole at time $s$. It is also easy to establish the inequality
$(M_s+1)(H_s+1)\geq D_s$. So from the result for $\max D_s$ in Theorem \ref{symmetricdistancetheorem},
it follows that $\liminf\PP(\max_{s\in[0,t]} M_s >\sqrt{\alpha t})$ is at least $1/2$,
and naturally one expects that this probability should in fact converge to 1. 
In the other direction, Arratia's result in (\ref{arratia}) together with Lemma \ref{dominationlemma}
would give an upper bound of $\sqrt{t\log t}$. 
\item We didn't determine an actual rate of linear growth for the maximal distance
of the bridge in the symmetric case. 
It seems unlikely that either 
the lower bound $\frac3{4\pi^2}$ or the upper bound $\frac14$
in Theorem \ref{symmetricdistancetheorem} are tight.
More generally, does there exist a deterministic asymptotic shape,
under appropriate rescaling, for the state
of the bridge at time $t/2$ (say)? 
That is, is there a constant limit in probability for the number of particles to the right
of $\beta\sqrt{t}$ at time $t/2$, divided by $\sqrt{t}$, for each $\beta$?
\item In the unconditioned process, 
for which $p$ and $t$ is it the case that $O$ is the most likely state at time $t$?
\end{enumerate} 

\section*{Acknowledgments}
Many thanks to Itai Benjamini for telling us about this problem 
(and for inventing the title).


\end{document}